\documentclass[a4paper,11pt,oneside]{article}

\usepackage{amsmath,amsthm,amsfonts,amssymb,amscd}
\usepackage{mathtools}
\usepackage{pgf,tikz}
\usepackage{float}
\usepackage{enumerate}
\usepackage{booktabs}
\usepackage{pgfplots}
\usepackage[noend]{algpseudocode}
\usepackage{graphicx,bm,xcolor}
\usepackage{algorithm}
\usepackage{algpseudocode}
\usepackage{hyperref}
\usepackage{caption}
\usepackage[T1]{fontenc}
\usepackage{t1enc}
\usepackage[polish,german,english]{babel}
\usepackage{verbatim} 

\usepackage{url}
\usepackage{authblk}

\usepackage[top=2cm, bottom=2cm, left=2cm, right=2cm]{geometry}

\theoremstyle{plain} 

\theoremstyle{definition} %

\theoremstyle{remark} %

\usepackage{type1cm}        
\usepackage{makeidx}         
\usepackage{graphicx}        
\usepackage{multicol}        
\usepackage[bottom]{footmisc}
\usepackage{algorithm}
\usepackage{algpseudocode}
\usepackage{tikz}
\usepackage{pgfplots}
\pgfplotsset{compat=1.18}

\usepackage{wrapfig}
\usepackage{newtxtext}       %
\usepackage[varvw]{newtxmath}       

\makeindex             

\newtheorem{Problem}{Problem}

\begin{document}

\title{Matrix-Free Geometric Multigrid Preconditioning Of Combined Newton-GMRES
For Solving Phase-Field Fracture With Local Mesh Refinement}
\author[1]{L. Kolditz}
\author[1]{T. Wick}

\affil[1]{Leibniz Universit\"at Hannover, Institut f\"ur Angewandte
  Mathematik, AG Wissenschaftliches Rechnen, Welfengarten 1, 30167 Hannover, Germany}

\date{}

\maketitle
	
\begin{abstract}
In this work, the matrix-free solution of quasi-static phase-field fracture problems
is further investigated. More specifically, we consider a quasi-monolithic formulation in which the irreversibility constraint is imposed with a primal-dual active set method. The resulting 
nonlinear problem is solved with a line-search assisted Newton method. Therein, the arising 
linear equation systems are solved with a generalized minimal residual method (GMRES), which is 
preconditioned with a matrix-free geometric multigrid method including geometric local mesh refinement. Our solver is substantiated with a numerical test on locally refined meshes.
\end{abstract}

\section{Introduction}
\label{sec:1}
This work is devoted to the efficient linear solution within the nonlinear 
solver of quasi-static phase-field fracture problems. Phase-field fracture 
remains a timely topic with numerous applications. Therein, vector-valued 
displacements and a scalar-valued phase-field variable couple. Moreover, the phase-field variable is subject to a crack irreversibility constraint. Due to nonlinear 
couplings, nonlinear constitutive laws and the previously mentioned 
irreversibility constraint, the overall coupled problem is nonlinear. Here, a 
line-search assisted Newton method is employed. Therein, the linear solution often 
is a point of concern. 

Based on prior work \cite{JoLaWi20}, a GMRES (generalized 
minimal residual method) iterative linear solver is employed. This is preconditioned 
with a matrix-free geometric multigrid (GMG) method. In the matrix-free context, the 
system matrix is not fully assembled \cite{KroKor2012}, 
which reduces the memory consumption. At the same time, this limits 
the choice of available smoothers for the multigrid preconditioner.
Here, a Chebyshev-Jacobi smoother is employed 
as it only requires matrix-vector 
products and an estimate of the largest eigenvalue. Moreover, the 
inverse diagonal of the system matrix 
is required, which however can be obtained from the local assembly without requiring the assembly of the entire matrix. With these ingredients 
a matrix-free solution can be set up.
Recent works of matrix-free solvers include 
problems in finite‐strain hyperelasticity \cite{DaPeAr20}, 
phase-field fracture  \cite{JoLaWi20_parallel,JoLaWi20},
Stokes \cite{jodlbauer2022matrixfree},
generalized Stokes \cite{WICHROWSKI2022101804},
fluid-structure interaction \cite{wichrowski2023exploiting}, 
discontinuous Galerkin \cite{KrKo19},
compressible Navier-Stokes equations \cite{guermond2021implementation}, 
incompressible Navier-Stokes and Stokes equations \cite{FrCaAnPa20} as well as 
sustainable open-source code developments \cite{MuKoKr20,CleHeiKaKr21}, 
matrix-free implementations on locally refined meshes \cite{MuHeiSaaKr22_arxiv},
implementations on graphics processors \cite{KroLjun19},
and performance-portable methods on CPUs and GPUs applied to solid mechanics \cite{Brown_etal:2022arXiv2204.01722}.
In the prior work \cite{JoLaWi20}, the solver only allowed globally refined meshes.
The main contribution of this work is that we combine a GMG preconditioner with locally refined 
meshes and primal-dual active set for the inequality constraint of phase-field fracture using local smoothing 
\cite{janssen2011adaptive}.
Our overall numerical solver is applied 
to one numerical test, namely 
Sneddon's example that is nowadays considered 
as a benchmark \cite{schroder2021selection}.
The outline of these conference proceedings 
is as follows. In Section \ref{sec:2}, the phase-field problem is stated.
Then, in Section \ref{sec_sol}, the numerical solutions are explained, namely 
the nonlinear solution via a combined Newton method and the linear solution 
with GMRES and matrix-free geometric multigrid preconditioning. Finally,
in Section \ref{sec_test} a numerical test substantiates our developements.

\section{Problem Formulation: The Phase-Field Model}\label{sec:2}
This section introduces the problem formulation for phase-field fracture which was originally formulated as an energy minimization problem \cite{francfort1998revisiting}. However, our starting point are the variational Euler-Lagrange equations. For this, we provide basic notations: given a sufficiently smooth material $\Omega\subset \mathbb{R}^d$, $d = 2$, the scalar-valued and vector-valued $L^2$-products over a smooth bounded domain $G\subset \Omega$ are defined by
\begin{equation}
    (x,y)_{L^2(G)} \coloneqq \int_{G} x\cdot y dG, \quad (X,Y)_{L^2(G)} \coloneqq \int_G X : Y dG,
\end{equation}
where $ X : Y$ denote the Frobenius product of two vector fields. If there is no subscript provided, the $L^2$-product over the whole material domain $\Omega$ is meant. Energy minimization problems in the phase-field fracture context (e.g. \cite{bourdin2000revisited,MiWheWi19})  usually consist of a displacement variable $u : \Omega \rightarrow \mathbb{R}^d$ and a phase-field variable $\varphi : \Omega \rightarrow [0,1]$. The phase-field variable $\varphi$ can be understood as an indicator function: It is defined such that $\varphi = 0$ where the material is fully broken and $\varphi = 1$ where it is intact. Only allowing a fully broken or a completely intact domain leads to discontinuities, which is treated by the Ambrosio-Tortorelli approximation. 
With this, we introduce a transition zone where $0<\varphi<1$ of width $2l$. In the following $l$ is called the length-scale parameter. Deriving the Euler-Langrange equations from computing directional derivatives, the solution sets on the 
continuous level are given by
        $\mathcal{V} \coloneqq H^1_0(\Omega),
        \mathcal{K}^n \coloneqq \{\psi \in \mathcal{W}\, \vert\, \psi - \varphi^{n-1} \leq 0 \text{ a.e. in } \Omega\}$, 
        where  $\mathcal{W} \coloneqq H^1(\Omega)$.
    The function space $\mathcal{K}^n$ is a convex set arising from the crack irreversibility constraint $\partial_t \varphi \leq 0$. In the quasi-static setting, this translates to $\varphi^n \leq \varphi^{n-1}$ with $\varphi^n:=\varphi(t_n)$ and $\varphi^{n-1}:=\varphi(t_{n-1})$. This results in an incremental grid $t_0,\dots t_N$ with the step-size $k_n = t_n - t_{n-1}$, where $t_0$ is the initial configuration and $t_N$ the end-time configuration.
The Euler-Lagrange equations are then given by \cite{MiWheWi19,wick2020multiphysics, KoWiMa2023modified}:
\begin{Problem}
For some given initial value $\varphi^0$ and for the incremental steps $t_n$ with $n=1,...,N$, find $(u^n,\varphi^n) \in \mathcal{V} \times \mathcal{K}^n$ such that
for all $\psi^u \in \mathcal{V}$ and $\psi^\varphi \in \mathcal{K}^n \cap L^\infty(\Omega)$
\begin{align}
\begin{split}
    \left(g(\varphi^{n})\sigma(u^n),e(\psi^u)\right) + ((\varphi^n)^2p^n, \operatorname{div} \psi^u) &= 0,
    \label{eq:u}
    \end{split}\\
    \begin{split}   
    (1-\kappa)(\varphi^n\sigma(u^n) : e(u^n), \psi^\varphi - \varphi^n) +2(\varphi^n p^n \operatorname{div} u^n, \psi^\varphi-\varphi^n) &\quad \\
    \qquad + G_C\left(\frac{1}{l}(1-\varphi^n,\psi^\varphi-\varphi^n) + l(\nabla\varphi^n,\nabla (\psi^\varphi-\varphi^n))\right) &\geq 0,
    \end{split}
\end{align}
where $p \in L^\infty(\Omega)$ is a given pressure, $G_C>0$ is the critical energy release rate and $\kappa>0$ is a regularization parameter. 
A study on how to find a proper setting for $\kappa$ and the length-scale parameter $l$ is given in \cite{kolditz2022relation}. The classical stress tensor of linearized elasticity $\sigma(u)$ and the symmetric strain tensor $e(u)$ are given by 
\begin{equation}
    \sigma(u) \coloneqq 2\mu e(u) + \lambda \operatorname{tr}(e(u))I, \quad e(u)\coloneqq \frac{1}{2}\left( \nabla u + \nabla u^T \right),
\end{equation}
with the Lam\'e parameters $\mu > 0$ and $\lambda$ with $3\lambda + 2\mu>0$ and the identity Matrix $I$. Lastly, the degradation function $g(\varphi^n)$ is defined by $g(\varphi^n) := (1-\kappa)(\varphi^n)^2+\kappa$.
\end{Problem}
To enhance the robustness of the nonlinear solution, we linearize the degradation function following \cite{heister2015primal,KoWiMa2023modified}. Therein, the phase-field $\varphi^n$ is replaced by the old incremental step solution or an extrapolation using previous incremental step solutions. Equation \eqref{eq:u} reads then 
\begin{equation}
\left(g(\Tilde{\varphi}^{n})\sigma(u^n),e(\psi^u)\right) + ((\Tilde{\varphi}^n)^2p^n, \operatorname{div} \psi^u) = 0.
\end{equation}
The second difficulty of the above problem is the fact that we have to deal with a constraint variational inequality system (CVIS). This inequality system can be equivalently formulated as an equality constraint system with an additional complementarity equation \cite{KoWiMa2023modified}:
\begin{Problem}
    For a given initial condition $\varphi^0$ and for the incremental steps $t_n$ with $n=1,...,N$, find $U^n =\{u^n,\varphi^n \}\in \mathcal{V} \times \mathcal{W}$ and $\lambda^n \in \mathcal{N}_+$ such that
\begin{align*}
    A(U^n)(\Psi) + (\lambda^n, \psi^\varphi) &= 0\quad \forall \Psi = \{\psi^u,\psi^{\varphi} \}\in \mathcal{V} \times \mathcal{W}\cap L^\infty,\\
    C(\varphi^n,\lambda^n) &= 0\quad \text{a.e.\ in }\Omega,
\end{align*}
with 
\begin{equation}
    \begin{split}
    A(U^n)(\Psi) &\coloneqq\left(g(\Tilde{\varphi}^n)\sigma^+(u^n),e(\psi^u)\right) +((\Tilde{\varphi}^n)^2p^n, \operatorname{div} \psi^u)\\
    &\quad  + (1-\kappa)(\varphi^n\sigma^+(u^n) : e(u^n), \psi^\varphi - \varphi^n) +2(\varphi^n p^n \operatorname{div} u^n, \psi^\varphi) \\
    &\quad + G_C\left(\frac{1}{l}(1-\varphi^n,\psi^\varphi-\varphi^n) + l(\nabla\varphi^n,\nabla (\psi^\varphi))\right),
    \end{split}
\end{equation}
and
\begin{equation}
    C(\varphi^n, \lambda^n)\coloneqq \lambda^n -\max\{0,\, \lambda^n +c(\varphi^n-\varphi^{n-1})\}.
\end{equation}
The solution set $\mathcal{N}_+$ for the Lagrange multiplier is defined by 
\begin{equation}
    \mathcal{N}_+\coloneqq \left\{\mu \in L^2(\Omega) \,|\,  ( \mu, v )_{L^2(\Omega)}\leq 0\quad \forall v \in L^2_-(\Omega) \right\}.
\end{equation}
\end{Problem}
This formulation is the starting point for a primal-dual active set (PDAS) method, which we use to treat the irreversibility condition. The idea is to split the domain, based on the structure of the complementarity condition, into two subdomains: The active set $\mathcal{A}$, where the constraint is active, i.e. the phase-field variable does not change, and the inactive set $\mathcal{I}$, where the constraint is inactive. In the latter, the problem can be treated and solved as an unconstrained problem. The active and inactive sets at each incremental step are defined such as 
\begin{align}\label{eq:AcInS}
    \lambda^n +c(\varphi^n-\varphi^{n-1}) \leq 0 \text{ a.e.\ in } \mathcal{I}^n, \quad
    \lambda^n +c(\varphi^n-\varphi^{n-1}) > 0 \text{ a.e.\ in } \mathcal{A}^n.
\end{align}
    The active set constant $c$ can be chosen arbitrarily. However, in other contributions \cite{hueber2005primaldual,popp2009finite,schroeder2016semismooth}, the authors state that it can have an influence on the performance. 
    Following our prior work \cite{KoWiMa2023modified}, we set $c$ as
   $
     c = c^k \coloneqq 2\left|\lambda_{h,i}^k / (\varphi_{h,i}^k-\varphi_{h,i}^{\operatorname{old}}) \right|.
$

\section{Nonlinear Solution With Inner Linear GMRES Iterative Solver And Matrix-Free Geometric Multigrid Preconditioning}
\label{sec_sol}
In this section, 
as discretization, we employ a finite element method with $H^1$-conforming bilinear finite elements for both the displacement and the phase-field, as an outer solver, we employ a combined Newton active set method and as an inner linear solver, we utilize a GMRES algorithm together with a geometric multigrid preconditioner. The newly developed code is implemented in a matrix-free framework using the finite element library \texttt{deal.II} \cite{dealii2021version94}. Related work 
having deal.II as well as basis was done in \cite{JoLaWi20,JoLaWi20_parallel}.

\subsection{A Combined Newton Type Algorithm}
In each incremental step, the Newton iteration to solve for $U^n \coloneqq \{u^n,\varphi^n\}\in \mathcal{V} \times \mathcal{W}$ and $\lambda^n \in \mathcal{N}_+$ is given by
\begin{equation*}
    A'(U^{n,k})(\delta U^{n,k+1}, \Psi) + (\lambda^{n,k},\psi^\varphi) = -A(U^{n,k})(\Psi) \quad \forall \Psi \in \mathcal{V}\times \mathcal{W},
\end{equation*}
which is solved with respect to 
\begin{equation*}
    C(\varphi^{n,k} + \delta \varphi^{n,k+1},\lambda^{n,k}) =0\quad \text{ a.e.\ in } \Omega,
\end{equation*}
for the Newton update $\delta U^{n,k+1}$ and the Lagrange multiplier $\lambda^{n,k+1}$. The solution is then updated via
\begin{equation*}
    U^{n,k+1} = U^{n,k} + \delta U^{n,k+1}.
\end{equation*}
The Jacobian $ A'(U^{n,k})(\delta U^{n,k+1}, \Phi)$ is given by
\begin{equation*}
    \begin{split}
    A'(U^{n,k}) &(\delta U^{n,k+1}, \Psi)
                = \left(g(\Tilde{\varphi}^{n}) \sigma(\delta u^{n,k+1}), e(\psi^u) \right) \\
                & + (1-\kappa)\left(\delta \varphi^{n,k+1} \sigma(u^{n,k}) : e(u^{n,k}) + 2 \varphi^{n,k} \sigma(\delta u^{n,k+1}) : e(u^{n,k}), \psi^\varphi\right) \\
                & + 2p\left(\delta\varphi^{n,k+1} \operatorname{div} u^{n,k} + \varphi^{n,k} \operatorname{div} \delta u^{n,k+1}, \psi^\varphi\right) \\
                & + G_C \left( \frac{1}{l}(\delta \varphi^{n,k+1}, \psi^\varphi) + l (\nabla \delta \varphi^{n,k+1}, \nabla \psi^\varphi)\right).
    \end{split}
\end{equation*}
In combination with an iteration on the active set, we obtain the scheme
outlined in Algorithm \ref{alg:active_set_cont}.
\begin{algorithm}
\caption{(Primal-dual active set method)}\label{alg:active_set_cont}
\begin{algorithmic}[1]
\State Set iteration index $k=0$
\While{$\mathcal{A}^{n,k} \neq \mathcal{A}^{n,k+1}$}
\State Determine the active set $\mathcal{A}^{n,k}$ and inactive set $\mathcal{I}^{n,k}$ with \eqref{eq:AcInS}
\State Find $\delta U^{n,k+1} \in \mathcal{V}\times \mathcal{W}$ and $\lambda^{n,k+1}\in \mathcal{N}_+$ with solving
\begin{align*}
    A'(U^{n,k})(\delta U^{n,k+1}, \Phi) + (\lambda^{n,k+1},\psi) &= -A(U^{n,k})(\Phi),\quad \forall \Phi\coloneqq \{v,\psi\} \in \mathcal{V}\times \mathcal{W},\\
    \delta \varphi^{n,k+1} &= 0\quad \text{ on } \mathcal{A}^k,\\
    \lambda^{n,k+1} &= 0\quad \text{ on } \mathcal{I}^k.
\end{align*}
\State Update the solution to obtain $U^{n,k+1}$ via
\begin{equation*}
    U^{n,k+1} = U^{n,k} + \delta U^{n,k+1}.
\end{equation*}
\State Update iteration index $k=k+1$
\EndWhile
\end{algorithmic}
\end{algorithm}
The combined Newton active set algorithm has two stopping criteria which have to be fulfilled: On one hand, the Newton residual has to be small enough while on the other hand the active set has to remain unchanged over two consecutive iterations. With the active set convergence criterion, we ensure that we indeed applied the constraints in the right way. Usually, the first prediction of the active set is very bad. Thus, we provide the final active set from the previous incremental step as initial active set for the next one.

\subsection{A Geometric Multigrid Block Preconditioner For GMRES}
In a matrix-free framework, only iterative methods which solely rely on matrix-vector products are applicable as smoothers. We choose a Chebyshev-accelerated polynomial Jacobi smoother. This method requires to precompute the inverse diagonal entries of the system matrix (i.e. the Jacobian) and an estimate for the eigenvalues. The latter can be obtained by employing a conjugate gradient method. On the levels where we are mainly interested in smoothing out the highly oscillating error parts (i.e. all levels $>0$), it suffices to compute the maximal eigenvalue and then approximate the smoothing range by $[\lambda_{\min}, \lambda_{\max}] \approx [0.08\lambda_{\max},1.2 \lambda_{\max}]$. 

On the coarsest grid (level $0$), where we solve the problem using the Chebyshev Jacobi smoother, we then compute both the maximal and the minimal eigenvalue for indeed solving the problem. We apply the preconditioner on the whole block system by performing one V-cycle on each of the symmetric positive definite diagonal blocks as follows. 
To neglect hanging node constraints during preconditioning, we employ a local smoothing approach, where the smoother only acts on the subdomains which are refined to the current level~\cite{janssen2011adaptive}.
We want to solve the following preconditioned inner linear system arising from finite element discretization and Newton's method
\begin{equation}
    P^{-1}\Tilde{G}\delta U=P^{-1}\Tilde{R}, \;
G \coloneqq 
    \begin{bmatrix}
        \Tilde{G}^{uu} & 0 \\
        \Tilde{G}^{u\varphi} & \Tilde{G}^{\varphi\varphi},
    \end{bmatrix},
    \;
P^{-1} \coloneqq 
    \begin{bmatrix}
        MG(\Tilde{G}^{uu}) & 0 \\
        0 & MG(\Tilde{G}^{\varphi\varphi})
    \end{bmatrix}
    ,
\end{equation}
with the mentioned Jacobi-Block smoother applied to each diagonal block.
    One challenge of the multigrid preconditioner is to transfer the information of possible constraints onto the coarser grids. In our case, we have to deal with three different types of constraints: boundary conditions, active set constraints and hanging node constraints. The latter are avoided due to local smoothing. The boundary constraints are transported to coarser meshes in a canonical manner. The biggest difficulty is to transfer the active set. Simply categorizing those dofs as active which are active on the next coarser level leads to too much information loss. Thus, we transfer the active set to each level following e.g. \cite{Hackbusch1983}: on level $k$, we only categorize a degree of freedom as inactive, if all its direkt neighbours are inactive, otherwise it is categorized as active.

\subsection{The Matrix-Free Approach And Final Algorithm}
The full Algorithm \ref{alg:full_algorithm} is realized in a matrix-free framework to reduce the memory consumption. From the implementation point of view, the concept is simple: Instead of assembling and storing the system matrix, we implement a linear operator which represents the application of the matrix to a vector. For this, the global matrix-vector-product can be split into local matrix-vector-products corresponding to the underlying finite elements. 
As previously mentioned, the implementation is realized with \texttt{deal.II}, which offers a toolbox of functionalities considering matrix-free finite-element approaches \cite{KroKor2012}.

\begin{algorithm}
\caption{full algorithm}\label{alg:full_algorithm}
\begin{algorithmic}[1]
\State Setup the system                                 \Comment{initialize grid $\mathcal{T}_h$, parameters, etc.}
\For{$n=1,2,...$}\Comment{timestep loop}
    \While{$\left(\mathcal{A}^{n,k-1} \neq \mathcal{A}^{n,k}\right)$ \textbf{or} $\left(\Tilde{R}(U_h^{n,k}) > \text{TOL}_N\right)$}
        \State Compute the active set $\mathcal{A}^{n,k} = \left\{ x_i  \bigg| \left[B\right]_{ii}^{-1} \left[R(U_h^{n,k})\right]_i +                   c(\varphi_{h,i}^{n,k} - \varphi_{h,i}^{n-1}) > 0\right\} $
        \State Set the active set constraints for the Newton update: $\delta \varphi_i^{n,k+1} = \varphi_i^{n-1}-\varphi_i^{n,k}$
        \State Set constraints for Newton update and assemble the residual $R(U_h^{n,k})$
        \State Setup the Multigrid preconditioner 
        \State Solve the inner linear system with GMRES and GMG preconditioner 
        \State Distribute the constraints on the Newton update 
        \State Choose maximum number of line search iterations $l_{\max}$
        \State Choose line search damping parameter $0 < \omega \leq 1 $
        \For{$l=1:l_{\max}$}
            \State Update the solution with $U_h^{n,k+1} = U_h^{n,k} + \delta U_h^{n,k+1}$
            \State Assemble the new residual $\tilde{R}(U_h^{n,k+1})$
        \If{$\lVert \tilde{R}(U_h^{n,k+1}) \rVert_2 < \lVert \tilde{R}(U_h^{n,k}) \rVert_2$}
            \State \textbf{break}
        \Else
            \State Adjust the Newton update with $\delta U_h^{n,k+1} := \omega^l\delta U_h^{n,k+1}$
        \EndIf
    \EndFor
    \State Update iteration index $k=k+1$
    \EndWhile
\EndFor
\end{algorithmic}
\end{algorithm}
\section{Numerical Test: Sneddon's Benchmark}
\label{sec_test}
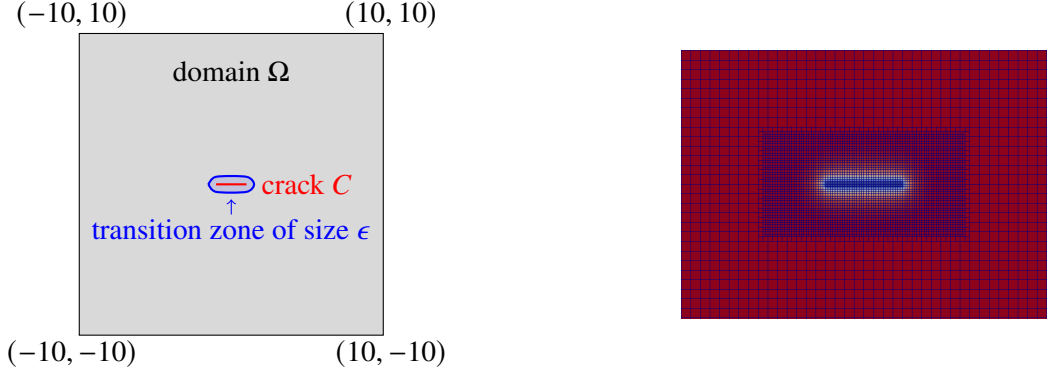
\begin{figure}[htbp!]
\begin{minipage}{0.47\textwidth}
\centering
\begin{tikzpicture}[scale = 0.9]
\node  at (-0.1,4.25) {$(-10,10)$};
\node  at (-0.1,-0.25) {$(-10,-10)$};
\node at (4.1,-0.25) {$(10,-10)$};
\node at (4.1,4.25) {$(10,10)$};
\draw[fill=gray!30] (0,0)  -- (0,4) -- (4,4)  -- (4,0) -- cycle;
\node at (2.0,3.5) {domain $\Omega$};
\draw [fill,opacity=0.1, blue] plot [smooth] coordinates { (1.7,2) (1.8,1.9) (2.2,1.9) (2.3,2) (2.2,2.1) (1.8,2.1) (1.7,2)};
\draw [thick, blue] plot [smooth] coordinates { (1.7,2) (1.8,1.9) (2.2,1.9) (2.3,2) (2.2,2.1) (1.8,2.1) (1.7,2)};
\draw[thick, draw=red] (1.8,2)--(2.2,2);
\node[red] at (3,2) {crack $C$};
\draw[->,blue] (2,1.6)--(2,1.8);
\node[blue] at (2,1.4) {transition zone of size $\epsilon$};
\end{tikzpicture}
\end{minipage}
\hspace{0.2cm}
\begin{minipage}{0.47\textwidth}
\centering
\includegraphics[scale = 0.25]{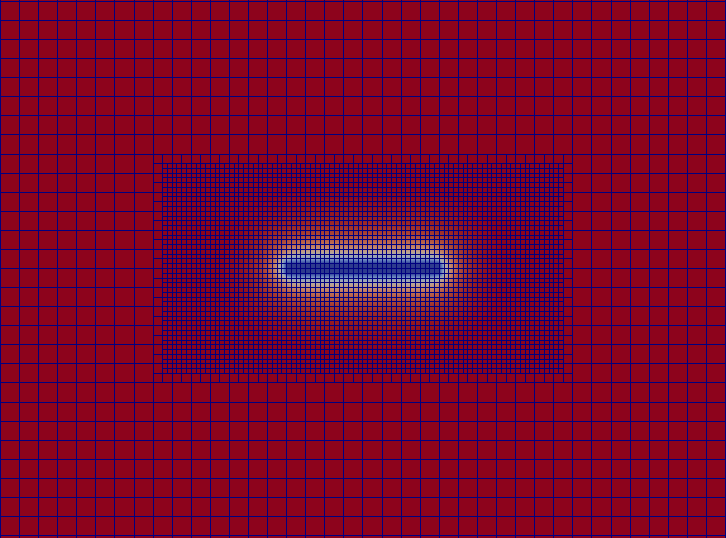}
\end{minipage}
    \caption{Left: Geometry two dimensional Sneddon test. Right: Geometrically locally refined mesh.
    }
    \label{geo_sneddon}
\end{figure}

In this section, we present numerical results for a stationary two dimensional benchmark test, where a one dimensional crack is prescribed in the center of the domain and a constant pressure is applied in the inner of the fracture. This test is also called the Sneddon benchmark test \cite{SneddLow69}. The two dimensional domain is given by $\Omega = (-10,10)^2$ as depicted in Figure \ref{geo_sneddon}. 

The fracture has a constant half length of $l_0=1.0$ and a varying width depending on the minimal element diameter. The parameters are given in Table \ref{tab:1}.

\begin{table}[!t]
\caption{The setting of the material and numerical parameters for the Sneddon 2d test.}
\label{tab:1}       
%
%
\centering
\begin{tabular}{p{2cm}p{5cm}p{2.5cm}}
\hline\noalign{\smallskip}
Parameter & Definition & Value \\
\noalign{\smallskip}\noalign{\smallskip}
$\Omega$ & Domain & $(-10,10)^2$  \\ 
$ h $ & Diagonal cell diameter & test-dependent\\
$ l_0 $ & Half crack length & $1.0$ \\ 
$ G_C $ & Material toughness & $1.0$ \\ 
$ E $ & Young's modulus & $1.0$ \\ 
$\mu$ & Lamé parameter &  $0.42$\\ 
$\lambda$ & Lamé parameter & $0.28$ \\ 
$ \nu $ & Poisson's ratio & $0.2$ \\ 
$ p $ & Applied pressure & $10^{-3}$ \\ 
$ l $ & length scale parameter & $2h$ \\
$ \kappa $ & Regularization parameter & $10^{-12}h$\\ 
& Number of global refinements & {$2$}\\
& Number of local refinements & {$0-8$}\\
{$\operatorname{TOL}_N$}& {Tolerance outer Newton solver} & {$10^{-7}$}  \\
{$\operatorname{TOL}_{LS}$}& {Tolerance GMRES} & {$\max{\{10^{-12},10^{-8}\Tilde{R}\}}$}  \\
\noalign{\smallskip}\hline\noalign{\smallskip}
\end{tabular}
\end{table}
The quantities of interest in this test case are given by the so-called total crack volume 
\[
\operatorname{TCV} \coloneqq \int_\Omega u(x,y) \nabla \varphi (x,y)\, d(x,y)
\]
and the crack opening displacement 
\[
\operatorname{COD}(x) \coloneqq [u \cdot n](x) \approx \int_{-\infty}^{\infty} u(x,y) \cdot \nabla \varphi(x,y) \, dy.
\]
The analytical solutions~\cite{SneddLow69} are given by $\operatorname{TCV}_{\operatorname{ref}} = \frac{2\pi \rho l_0^2}{E'}$ and $\operatorname{COD}_{\operatorname{ref}} = 2\frac{pl}{E'}\left(1-\frac{x^2}{l_0^2} \right)^\frac{1}{2}$.

\begin{table}
\caption{Computational results for the Sneddon test with $2$ global and $0-8$ local refinements.}
\label{tab:2}       
%
%
\centering
\begin{tabular}{p{1cm}p{1.3cm}p{1.5cm}p{2.0cm}p{1.5cm}p{1.8cm}p{1.8cm}}
\hline\noalign{\smallskip}
$h$ &$\#$DoFs  & TCV error&  $\#$Newton steps &$\varnothing$ lin. iter. & $\max$ $\#$lin. iter.&Wall-time [s]\\ 
$0.7071$&$5043$&$288.00\%$&$19$&$1.95$&$6$&1.20\\
$0.3536$&$5745$&$115.00\%$&$21$&$1.71$&$7$&1.58\\
$0.1768$&$8445$&$45.90\%$&$20$&$1.80$&$8$&2.29\\
$0.0884$&$16953$&$18.70\%$&$24$&$1.70$&$9$&4.31\\
$0.0442$&$48321$&$7.67\%$&$25$&$2.08$&$10$&9.84\\
$0.0221$&$168609$&$2.88\%$&$29$&$2.69$&$13$&36.20\\
$0.0110$&$639537$&$0.62\%$&$29$&$3.21$&$18$&132.00\\
$0.0055$&$2502945$&$0.50\%$&$33$&$3.88$&$25$&633.00\\
$0.0028$&$9916113$&$1.00\%$&$30$&$5.40$&$40$&2800.00\\

\noalign{\smallskip}\hline\noalign{\smallskip}

\noalign{\smallskip}\hline\noalign{\smallskip}
\end{tabular}
\end{table}
\begin{figure}[!t]
    \centering
    \includegraphics[width=0.98\linewidth]{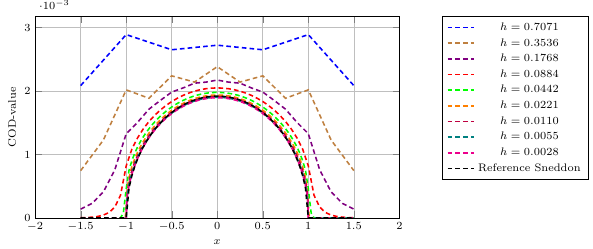}
    \caption{Visualization of the $\operatorname{COD}$-values for different $h$. The corresponding exact $\operatorname{TCV}$ values are given in Table \ref{tab:2}.}
    \label{fig:cod}
\end{figure}
In Table \ref{tab:2} and Figure \ref{fig:cod}, the results of the TCV and the COD on nine different refinement levels are depicted and compared to the reference solution. All the computations were performed on the same machine and on 4 cores. We can clearly observe that the numerical solution tends to the reference solution under grid refinement. The number of linear iterations indicates robustness under  mesh refinement. However, there still appear peaks, which assumingly come up due to the bad prediction of the active set in the initial Newton step. Once the active set is well predicted, the number of linear iterations is very stable under (local) mesh refinement.

\bibliographystyle{abbrv}

\end{document}